\documentclass[11pt]{article}
\usepackage{amssymb,amsmath}
\usepackage{graphics}
\usepackage{float}
\usepackage{subfigure}
\usepackage{color,fullpage}
\usepackage{bm}
\usepackage[title,titletoc]{appendix}
\usepackage{algorithm,algorithmic}  
\usepackage{multirow}
\usepackage{pgfplots}
\usepackage{booktabs}
\usepackage{amsthm}
\usepackage{mathrsfs}

\usepackage{flushend} 
\usepackage{ulem}    
\usepackage[numbers, sort]{natbib}
\usepackage{enumerate}

\newtheorem{assumption}{Assumption}[section]

\newtheorem{theorem}{Theorem}[section]
\newtheorem{lemma}{Lemma}[section]
\newtheorem{definition}{Definition}[section]
\newtheorem{example}{Example}[section]

\newtheorem{remark}{Remark}[section]

\makeatletter

\@addtoreset{equation}{section} \makeatother

\begin{document}

\title{Nonlinear preconditioned primal-dual method for a class of structured minimax problems}

\author{
Lu Zhang\thanks{Department of Mathematics, National University of Defense Technology,
	Changsha, Hunan 410073, China. Email: \texttt{zhanglu21@nudt.edu.cn}}
\and Hongxia Wang \thanks{
Department of Mathematics, National University of Defense Technology,
Changsha, Hunan 410073, China.  Email: \texttt{wanghongxia@nudt.edu.cn}
}
\and Hui Zhang\thanks{Corresponding author.
Department of Mathematics, National University of Defense Technology,
Changsha, Hunan 410073, China.  Email: \texttt{h.zhang1984@163.com}
}
}


\date{\today}

\maketitle

\begin{abstract}
	We propose and analyze a general framework called nonlinear preconditioned primal-dual with projection for solving nonconvex-nonconcave and non-smooth saddle-point problems. The framework consists of two steps. The first is a nonlinear preconditioned map followed by a relaxed projection onto the separating hyperspace we construct. One key to the method is the selection of preconditioned operators, which tailors to the structure of the saddle-point problem and is allowed to be nonlinear and asymmetric. The other is the construction of separating hyperspace, which guarantees fast convergence. This framework paves the way for constructing nonlinear preconditioned primal-dual algorithms. We show that weak convergence, and so is sublinear convergence under the assumption of the convexity of saddle-point problems and linear convergence under a metric subregularity.
	We also show that many existing primal-daul methods, such as the generalized primal-dual algorithm method, are special cases of relaxed preconditioned primal-dual with projection.
\end{abstract}

\section{Introduction}
Consider the nonconvex-nonconcave and non-smooth saddle-point problem
\begin{equation}
	\label{eq1:1}
	\min_{x\in \mathcal{X}}\max_{y\in \mathcal{Y}} \Phi (x,y):= f(x)+\phi (x,y)-g(y),
\end{equation}
where $f:\mathcal{X}\rightarrow\mathcal{H},g:\mathcal{Y}\rightarrow\mathcal{H}$ are proper, l.s.c., convex but not necessarily smooth functions, $\phi:\mathcal{H}\rightarrow\mathcal{H}$ is continuously differentiable, weakly-convex in $x$ and weakly-concave in $y$, and $\mathcal{X},\mathcal{Y}$ are general real Hilbert spaces (unconstrained case) or convex sets (constrained case).
Most of variational models arising in image processing problems \cite{chambolle2011first,chambolle2016ergodic,esser2010general}, Nash equilibrium problems in game theory \cite{von2009optimization,krawczyk2000relaxation}, and canonical convex programming problems \cite{he2020optimal} can be modeled as a saddle-point problem (\ref{eq1:1}). 

The primal-dual method is an elementary method initially designed for solving the bilinear 	case of (\ref{eq1:1})
\begin{equation}
\label{1:2}
\min_{x\in \mathcal{X}}\max_{y\in \mathcal{Y}} \Phi (x,y):= f(x)+\langle Kx,y\rangle-g(y),
\end{equation}
where $K$ is a linear operator with norm $$\|K\|=\max\{\|Ky\|:y\in\mathcal{Y},\|y\|\leq 1\}.$$
It is initially proposed in \cite{zhu2008} and referred to as primal-dual hybrid gradient (PDHG).
Drawing upon a linear extraplation based on the current and previous iterates, \cite{chambolle2011first} formally proposes the primal-dual method and it works by alternating dual and primal variables of the form 
\begin{equation}
\label{1:3}
\left\{
\begin{array}{l}
	x_{k+1}=\arg\min_x \Phi (x,y_k)+\frac{1}{2\tau}\|x-x_k\|^2\\
	\bar{x}_{k+1}=x_{k+1}+\theta(x_{k+1}-x_k)\\
	y_{k+1}=\arg\max_y \Phi(\bar{x}_{k+1},y)-\frac{1}{2\sigma}\|y-y_k\|^2
\end{array},
\right.
\end{equation}
where $\tau>0$ and $\sigma>0$ are regularization parameters, $\theta\in [0,1]$ is an extrapolation parameter. When $\theta=0$, the primal-dual method reduces to the Arrow-Hurwicz algorithm \cite{1960study}.
For problem (\ref{1:2}) in convex-concave setting (both $f$ and $g$ are convex),
\cite{chambolle2011first} proves the primal-dual method converges with rate $O(\frac{1}{N})$ for $\theta=1$ as long as $\tau$ and $\sigma$ are chosen such that $\tau\sigma \|K\|^2<1$.
An important observation in \cite{pock2011diagonal} is that for problems where the operator $K$ has a complicated structure, $\|K\|$ is hard to estimate or it might be so large that it slows down the convergence rate of the algorithm. 
To alleviate the shortcoming, we study primal-dual method for problem (\ref{eq1:1}) from the perspective of preconditioned methods in this paper. Drawing upon the preconditioned technique, we can avoid approximating the norm of operator $K$ and accelerate the convergence rate without changing the computational complexity of the iterations. We refer to, e.g., \cite{Bregman1967,lu2023unified,he2012convergence,MALITSKY2018,he2022generalized,bui2020warped} for the other perspectives of the primal-dual methods.

For problem (\ref{1:2}), there have been developed intensive primal-dual methods with different preconditioners. In particular, we can view the classic primal-dual method \cite{chambolle2011first} as a proximal point algorithm (PPA) with the preconditioner
\begin{equation}
\label{1:4}
M =
\left[
\begin{array}{cc}
\frac{1}{\tau}I & -K^T \\
-\theta K & \frac{1}{\sigma}I
\end{array}
\right]
.
\end{equation}
We observed that when $\theta=1$, the preconditioner $M$ is linear and symmetric; when $\theta\in [0,1)$, the preconditioner $M$ is linear and 
asymmetric. 
Note that the authors in \cite{chambolle2011first} did not provide the convergence results for the asymmetric case. The asymmetric structure of preconditioner makes it hard to guarantee convergence. As studied in various literature such as \cite{he2012convergence,he2022generalized}, convergence of the primal-dual method with the asymmetric preconitioner can be established only when the primal-dual method is corrected by a correction step with some restrictive  conditions on parameters. The correction steps in \cite{he2012convergence,he2022generalized} brings an additional advantage, namely relaxing the involved parameters. We refer the reader to \cite{pock2011diagonal,lu2023unified,liu2021acceleration,guo2023preconditioned} for more preconditioned primal-dual algorithms in bilinear case.

For problem (\ref{eq1:1}), 
existing preconditioners are also classified as symmetric and asymmetric. The paper \cite{bui2020warped} propose the preconditioned primal-dual with warped resolvent but does not give the concrete form of the preconditioner; the assumptions are not based on the special structure of the saddle-point problem (\ref{eq1:1}). 
The semi-anchored multi-step gradient descent algorithm (SA-MGDA) in \cite{lee2021semi}
 provides a nonlinear and symmetric preconditioner $M=\nabla h$ with 
\begin{eqnarray}
\label{1:6}
h(x,y)=\frac{1}{\tau}(\psi_x(x)+\psi_y(y))-\phi(x,y),
\end{eqnarray}
where $\tau>0$ is a relaxed parameter. \cite{lee2021semi} only provides the convergence result of iterations when $\phi$ is strongly convex-concave ($\phi(x,y)$ is strongly convex about $x$ and strongly-concave about $y$).
The generalized primal–dual
proximal splitting (GPDPS) method in \cite{Clason2021} provides a nonlinear and asymmetric preconditioner, which is a conceptually straightforward extension of the primal–dual method to saddle-point problems. The convergence result of GPDPS depend on some specific assumptions about the first and second derivative of $\phi$. In detail, the weak convergence requires the first derivative of $\phi$ to be maximally momotone; the strong convergence resuires the first derivative of $\phi$ to be strongly monotone or the second derivative of $\phi$ have better properties, e.g, second continuously differentiable, through Assumption 3.2 in \cite{Clason2021}. A natural question is whether we can design a new framework with a class of nonlinear and asymmetric preconditioners so that we can obtain convergence results under weaker assumptions than the aforementioned results and recover other existing algorithms?

Motivated by a class of preconditioned primal-dual algorithms in \cite{bui2020warped}, we propose a more general form of preconditioner $M$ for solving nonconvex-nonconcave saddle-point problems (\ref{eq1:1}), provided that certain conditions on parameters are required.
However, the nonlinear and asymmetric structure of preconditioner may break convergence and even lead to the implicit update of the algorithm.
Consequently, based on the separate and project principle in \cite{combettes2001fejer}, we strictly construct a separating hyperspace between the current iterate and the solution set and project the current iterate onto the separating hyperspace to ensure convergence. For the implicit update, on the one hand, we propose to use Bregman approximation \cite{Zhangh2022} to solve the possible implicit update; on the other hand, we directly construct the appropriate form of preconditioner to avoid the implicit update. For simplicity, we mainly use the second.
In this paper, different from other existing works, we place at the core of our study as constructing the concrete form of preconditioner and the separating hyperplane according to the special structure of the saddle-point problem (\ref{eq1:1}).

\textbf{Contributions.} The main contributions of this paper can be summarized as follows.
\begin{itemize}
	\item We propose a framework for the problem (\ref{eq1:1}), the nonlinear preconditioned primal-dual with projection method. To the best of knowledge,  taking the coupling of the primal variable and the dual variable into account, we present the most general preconditioner with nonlinear and asymmetric structure.
	Our framework is an implementation of \cite{bui2020warped,lee2021semi} with a concrete and more general form of the preconditioner $M$.
	\item To address the challenge of nonlinear and asymmetric preconditioners on convergence, we construct a correction step drawing upon the separate and project principle. This becomes a analytic framework of convergence of primal-dual with different preconditioners.
	\item Our framework can recover most of primal-dual algorithms under weaker assumptions. The interesting thing is that we also can generate new algorithm from our framework according to our needs.
\end{itemize}

\textbf{Organization.}
The paper is organized as follows. Section 2 is dedicated to notation and discussion about the selection of preconditioner $M$. In Section 3, we formally propose the nonlinear preconditioned primal-dual with projection method. The requirements of operators and parameters in the preconditioner are also discussed. Weak, linear, and sublinear convergence of the nonlinear preconditioned primal-dual with projection method are respectively analyzed under different conditions in Section 4. The connections to existing algorithms are discussed in Section 5. Section 6 is our conclusion.

\section{Preliminaries}
In this section, we summarize some knowledge of operators, which can refer to \cite{bauschke2011convex}. In particular, we reformulate the saddle-point problem (\ref{eq1:1}) into the operator inclusion problem and variational inequality (VI).

\subsection{The operators}
Throughout the paper, $\mathcal{X}$ and $\mathcal{Y}$ are general real Hilbert spaces $\mathcal{H}$.
The powerset of $\mathcal{H}$ is denoted $2^{\mathcal{H}}$. 
Let $P:\mathcal{H} \rightarrow 2^{\mathcal{H}}$ be a set-valued map.
The graph of $P$ is $\operatorname{gra}P=\{(x, u): u \in P x\}$. 
We denote by $\operatorname{dom} P=\{u\in\mathcal{H} : Pu\neq\emptyset\}$ the domain of $P$, by $\operatorname{ran} P=\{v\in\mathcal{H}: v\in Au,\forall u\in \mathcal{H} \}$ the range of $P$, by $\operatorname{zer} P=\{u\in\mathcal{H}:0 \in Pu\}$ the set of zeros of $P$, and by $P^{-1}$ the inverse of $P$, i.e. $\operatorname{gra} P^{-1}=\{(v,u):v\in Pu\}$. The inverse $P^{-1}$ is single-valued from $2^{\mathcal{H}}$ to $\mathcal{H}$ if $P$ is injective.
An operator $P$ is monotone if 
$$\langle u-v,Pu-Pv\rangle\geq 0,\forall u,v\in\mathcal{H},$$ 
and maximally monotone if, additionally, there exists no monotone operator $M:\mathcal{H} \rightarrow 2^{\mathcal{H}}$ such that $\operatorname{gra} P\subset \operatorname{gra} M$.  
An operator $P$ is self-adjoint if $\langle Pu,v\rangle = \langle u,Pv\rangle$ for all $u,v\in\mathcal{H}$. The projection operator onto the convex set $S$ is defined as 
\begin{equation*}
\Pi_{S}(x):=\arg\min_{y\in S}\{\|x-y\|\}.
\end{equation*}
Let $f$ be a proper, l.s.c. and convex function. The standard proximal resolvent of $f$ is defined as $(I+\partial f)^{-1}$ and corresponding proximal mapping is defined as 
$$
\text{prox}_f(y):=\arg\min_{x\in\mathcal{X}}\{f(x)+\frac{1}{2}\|x-y\|^2\}.
$$

\subsection{The reformulation of (\ref{eq1:1})}
Denote the primal-dual solution pair as $u=(x,y)$.
In constrained case, the problem (\ref{eq1:1}) can be reformulated into the following variantial inequality (VI): find $u^*=(x^*,y^*)\in \mathcal{X}\times \mathcal{Y}$ such that
\begin{equation}
\label{eq1:4}
l(u)-l(u^*)+(u-u^*)^TBu^*\geq 0,\forall u\in \mathcal{X}\times \mathcal{Y},
\end{equation}
where 
$$
l(u)=f(x)+g(y).
$$
For more details the reader is referred to \cite{he2022generalized}.

In the unconstrained case, finding the saddle point of problem (\ref{eq1:1}) is equivalent to finding a zero of the following set-valued operator $P$, namely,  
\begin{equation}
\label{eq1:2}
0\in Pu:=Au+Bu,u\in \mathcal{H}\times \mathcal{H},
\end{equation}
where
\begin{equation}
\label{eq1:3}
Au=
\left[
\begin{array}{c}
\partial f(x)\\
\partial g(y)
\end{array}
\right],
Bu=
\left[
\begin{array}{c}
\nabla_x\phi (u) \\
-\nabla_y\phi(u)
\end{array}
\right].
\end{equation}
Let $\operatorname{zer} P=\{u^* :0\in Pu^*\}$ be a nonempty solution set.
If $\phi$ is convex-concave, basic results in \cite{bauschke2011convex} show that $0\in Pu^*$ iff $u^*$ is a saddle point; if $\phi$ is continuously differentiable, then \cite{clarke1990optimization} indicates that $0\in Pu^*$ is a first-order necessary optimality condition for a saddle point $u^*$. Our primary aim of this study is to find a saddle point $u^*=(x^*,y^*)$ that satisfies the necessary optimality conditions of the problem (\ref{eq1:1}).

\section{The proposed preconditioner $M$}

We review that the warped resolvent of $P:\mathcal{H} \rightarrow 2^{\mathcal{H}}$ with respect to the preconditioner $M:\mathcal{H} \rightarrow\mathcal{H}$ proposed by \cite{bui2020warped} is defined as 
\begin{equation*}
R=(M+P)^{-1}M.
\end{equation*}
This reduces to the proximal resolvent operator $(I+P)^{-1}$ for $M=I$, where $I$ is indentity operator.
In this section, we explore a class of preconditioner $M$ with general structure and analyze the conditions that make the corresponding resolvent $R$ well-defined. We also provide several examples for our proposed preconditioner.
In order to make up for the shortcomings of \cite{lee2021semi} and couple primal variable $x$ and dual variable $y$, we introduce new operators into the preconditioner $M$. It will make the scheme fully implicit when new operators on both the lower left side and the upper right position in $M$. Actually, this is not feasible and makes it as difficult as the initial problem.
If we put the operator in the upper right side of the preconditioner, it is easy to see that the scheme can be easily made semi-implicit by exchanging the order between $x$ and $y$. 
Since primal variable $x$ and dual variable $y$ have the same status, there is no difference adding the operator in the upper right side or the lower left side.
Therefore, we simply put the operator $Q$ in the lower left side of the preconditioner.
We also consider generalizing $M$ in theory by employing operators $N_1, N_2$ instead of gradients generated from some functions. 
Thus, we propose two more general preconditioners $M$ based on the special structure of the saddle-point problem, namely,
\begin{equation}
	\label{3:1}
	M =
	\left[
	\begin{array}{cc}
	N_1 & 0 \\
	Q & N_2
	\end{array}
	\right]
	-
	\left[
	\begin{array}{c}
	\nabla_x\phi\\
	-\nabla_y\phi
	\end{array}
	\right],
\end{equation}
and
\begin{equation}
\label{3:2}
M =
\left[
\begin{array}{cc}
N_1 & 0 \\
Q & N_2
\end{array}
\right]
-
\left[
\begin{array}{c}
\nabla_x\phi\\
\nabla_y\phi
\end{array}
\right],
\end{equation}
where $N_1,N_2,Q$ represent operators rather than matrices.

\subsection{The well-defined resolvent}

Next, we turn to explore requirements to guarantee the single-valuedness of the warped resolvent $R$ with our proposed preconditioners. We refer to \cite{bui2020warped} for sufficient conditions that guarantee the warped resolvent are well-defined.
\begin{lemma}[\cite{bui2020warped}]
	For mappings $P:\mathcal{H}\rightarrow2^{\mathcal{H}}$ and $M:\mathcal{H}\rightarrow\mathcal{H}$, if $\operatorname{ran} M\subset \operatorname{ran}(P+M)$ and $P+M$ is injective, then $(M+P)^{-1}M$ is single-valued.
\end{lemma}

\begin{assumption}
	\label{assum:1}
	$f:\mathcal{X}\rightarrow\mathcal{H},g:\mathcal{Y}\rightarrow\mathcal{H}$ are proper, l.s.c., convex. The function $\phi:\mathcal{H}\rightarrow\mathcal{H}$ is continuously differentiable and has a Lipschitz continuous gradient. Further, $\phi(\cdot,y)$ is $L_{xx}$-weakly convex for fixed $y$, and $-\phi(x,\cdot)$ is $L_{yy}$-weakly convex for fixed $x$.
\end{assumption}


\begin{lemma}
	\label{lemma:2}
	Under Assumption \ref{assum:1}, $B$ and $P$ are $L$-Lipschitz continuous and $\gamma$-weakly monotone.
	
\end{lemma}

Throughout the paper, all proofs can be found in Appendix. Now that we have proposed assumptions on the operator $P$, we turn to consider assumptions about the preconditioner $M$.

\begin{assumption}
	\label{assum:2}
	$N_1$ is $\mu_1$-strongly monotone and $L_1$-Lipschitz continuous, $N_2$ is $\mu_2$-strongly monotone and $L_2$-Lipschitz continuous, and $N_1,N_2,Q$ are self-adjoint and single-valued. Further, $N_1,N_2$ are invertible. $M+P$ is maximally monotone. 
\end{assumption}

\begin{lemma}
	\label{lemma:3}
	Under Assumption \ref{assum:2}, $M$ is $\mu$-strongly monotone and $q$-Lipschitz continuous, where $\mu=\min\{\mu_1,\mu_2\}-L-\frac{\|Q\|}{2}>0$ and $q=\gamma+\max\{L_1,L_2\}+\frac{\|Q\|}{2}$.
\end{lemma}

Assumption \ref{assum:1} and Assumption \ref{assum:2} are sufficient for the guarantee that the resolvent $R$ is well defined, and we summarize as follows. 
\begin{lemma}
	\label{lemma:4}
	Suppose that Assumption \ref{assum:1} and Assumption \ref{assum:2} hold and $\mu>\gamma$, then
	\begin{enumerate}[(i)]
		\item $R$ is single-valued.
		\item $R$ has full domain.
		\item $\operatorname{fix} R=\operatorname{zer} P$,
		where $\operatorname{fix} R=\{x:x=Rx\},\operatorname{zer}P=\{x:0\in Px\}$.
	\end{enumerate}
\end{lemma}

The following analysis is based on Assumption \ref{assum:1} and Assumption \ref{assum:2} with $\mu>\gamma$. Given that $R$ has been well-defined, the preconditioned primal-dual method iteratively applies the preconditioner $M$ as 
\begin{equation}
\label{3:3}
u_{k+1}=Ru_k.
\end{equation} 
Equipped with the first preconditioner (\ref{3:1}), we can easily verify that the iterate generated by (\ref{3:3}) can be characterized as follows: 
\begin{equation}
\label{3:4}
\left\{
\begin{array}{l}
x_{k+1}=(\partial f+N_1)^{-1}(N_1 x_{k}-\nabla_x\phi(x_k,y_k))\\
y_{k+1}=(\partial g+N_2)^{-1}(N_2 y_k+\nabla_y\phi(x_k,y_k)+Q(x_{k}-x_{k+1}))
\end{array}.
\right.
\end{equation}
It follows from the simple lower-triangular structure of operator $M+P$ that preconditioned primal-dual method with preconditioner (\ref{3:1}) has an explicit update. 

With analogous reasoning, the preconditioned primal-dual with the second preconditioner (\ref{3:2}) iterates as:
\begin{equation}
\label{3:5}
\left\{
\begin{array}{l}
x_{k+1}=(\partial f+N_1)^{-1}(N_1 x_{k}-\nabla_x\phi(x_k,y_k))\\
y_{k+1}=(\partial g+N_2)^{-1}(N_2 y_k-\nabla_y\phi(x_k,y_k)+Q(x_{k}-x_{k+1})+2\nabla_y\phi(x_{k+1},y_{k+1}))
\end{array}.
\right.
\end{equation}
Facing with the challenge of implicit update about dual variable $y$ in (\ref{3:5}), we refer to two different approaches in \cite{lee2021semi,valkonen2021first}. Here, we also propose our own approach by weakening the requirement for  $\phi$. Assume that $(\phi,h)$ is Lipschitz-like convex \cite{Bauschke2017,Zhangh2022} and $h$ is a Legendre function. We use the Bregman approximation about $\phi(x_{k+1},y)$, i.e.,
\begin{equation}
\label{eq2:13}
\phi(x_{k+1},y)\approx \phi(x_{k+1},y_k)+\langle \nabla_y \phi(x_{k+1},y_k),y-y_k\rangle +\delta D_h(y,y_k),
\end{equation}
where $\delta>0$ is a regularized parameter and the Bregman distance associated to $h$ is defined by
$$D_h(x,y)=h(x)-h(y)-\langle \nabla h(y),x-y\rangle.$$
Then, the corresponding update about $y$ changes to
\begin{equation}
\label{eq2:14}
\begin{aligned}
y_{k+1}&=(\partial g+N_2-2\delta\nabla h)^{-1}((N_2-2\delta\nabla h)(y_k)-\nabla_y\phi(u_k)+2\nabla_y\phi(x_{k+1},y_k)+Q(x_k-x_{k+1})).
\end{aligned}
\end{equation}
When $h(x)=\frac{1}{2}\|x\|^2$, (\ref{eq2:13}) is the linearized approximation, and (\ref{eq2:14}) reduces to
$$
y_{k+1}=(\partial g+N_2-2\delta I)^{-1}((N_2-2\delta I)(y_k)-\nabla_y\phi(u_k)+2\nabla_y\phi(x_{k+1},y_k)+Q(x_k-x_{k+1})),
$$
which recovers the update about $y$ of SA-MGDA \cite{lee2021semi}.
In particular, the correction we propose about $y$-update in (\ref{eq2:14}) can relax the requirement for $\phi$ from Lipschitz gradient continuity to Lipschitz-like convexity.

\begin{remark}
	The $y$-update in (\ref{eq2:14}) has advantages over (\ref{3:4}) in terms of the weak requirement on $\phi$. In this paper, we will only explore the relatively simple preconditioner $M$ (\ref{3:1}),
	and leave the preconditioner $M$ (\ref{3:2}) for the future.
\end{remark}



The following aim of this paper is to ensure that the iterations generated by (\ref{3:4}) converge to the zero of the weakly monotone operator $P$.

\section{Nonlinear preconditioned primal-dual method}
In this section, we develop a nonlinear preconditioned primal-dual with projection method with the preconditioner $M$ (\ref{3:1}) for solving saddle-point problem (\ref{eq1:1}). We first elucidate the details of algorithmic design and then formally propose the algorithm in Subsection \ref{section4:1}. Subsection \ref{section4.2} provides examples of preconditioners.
Finally, the selection of parameters are discussed in Subsection \ref{section4.3}.
\subsection{The projected correction step}
\label{section4:1}

It is difficult to verify whether the preconditioned primal-dual method with the nonlinear and asymmetric preconditioner $M$ converges to the solution set $\operatorname{zer}P$. 
The challenge is that we have proposed a preconditioner $M$ with a versatile structure but we cannot guarantee convergence of the corresponding algorithm. 
Motivated by the separate and project principle and the geometric construction of F\'ejer monotonicity \cite{combettes2001fejer},
projections onto the solution set $\text{zer} P$ can be difficult. Hence, it is preferable to project onto a suitable separating halfspace, which contains the solution set $\text{zer} P$ but not the estimated iterate.

For all $z\in \operatorname{zer}P, u\in\mathcal{H}$.
Let $ C=M+B$.
From $0\in (A+B)z, Ru=(P+M)^{-1}Mu$, we have
$$
\begin{aligned}
0&\geq -\langle Az-ARu,z-Ru\rangle\\
&=-\langle -Bz+BRu+MRu-Mu,z-Ru\rangle\\
&=\langle Cu-CRu,z-Ru\rangle+\langle Bz-Bu,z-Ru\rangle\\
&=\langle Cu-CRu,z-Ru\rangle+\langle Bz-Bu,z-u\rangle+\langle Bz-Bu,u-Ru\rangle\\
&\geq \langle Cu-CRu,z-Ru\rangle+\frac{1}{L}\|Bz-Bu\|^2-\frac{1}{2\varepsilon}\|Bz-Bu\|^2-\frac{\varepsilon}{2}\|u-Ru\|^2\\
&=\langle Cu-CRu,z-Ru\rangle-\frac{L}{4}\|u-Ru\|^2,
\end{aligned}
$$
where we have used the monotonicity of operator $A$, the $\frac{1}{L}$-coercivity of $B$, Young's inequality, and $\varepsilon=\frac{L}{2}>0$. Thus, we have
$$
\langle (M+B)u-(M+B)Ru,z-Ru\rangle\leq  \frac{L}{4}\|u-Ru\|^2,\forall u\in\mathcal{H}.
$$

For every $u\in\mathcal{H}$, the separating halfspace can be redefined as
$$H_k=\{z:\psi_{u_k}(z)\leq 0\},$$
where
$$\psi_u(z)=\langle z-Ru,(M+B)u-(M+B)Ru\rangle-\frac{L}{4}\|u-Ru\|^2$$
and $L$ is the Lipschitz continuous constant of $P$ in Assumption \ref{assum:2}. Next, we show that the halfspace sequence $\{H_k\}_{k\in\mathbb{N}}$ can strictly separate the iterate sequence $\{u_k\}_{k\in\mathbb{N}} $ and the solution set $\operatorname{zer} P$ in each iterate.

\begin{theorem}
	\label{th1}
	Suppose that Assumption \ref{assum:1} and Assumption \ref{assum:2} hold, we have the following results about $\psi_u(z)$:
	\begin{enumerate}[(i)]
		\item $\psi_u(u)\geq (\mu-\gamma-\frac{L}{4})\|Ru-u\|^2\geq 0$ when $\mu>\gamma+\frac{L}{4}$, for all $ u\in \mathcal{H}$.
		\item $\psi_u(u)=0$ iff $u\in \operatorname{zer}P=\operatorname{fix} R$.
		\item $\psi_u(z)\leq 0$ for all $z\in \operatorname{zer}P=\operatorname{fix} R$ and $u\in\mathcal{H}$.
	\end{enumerate}
	
\end{theorem}
Theorem \ref{th1} shows that the separation is strict unless $u_k$ has already solved the operator inclusion problem. 
We have constructed a hyperspace separating the solution set and iterate, and then a relaxed projection onto this separating hyperplane finishes one iteration.
According to \cite{combettes2001fejer}, separation and projection are key tools that allow our preconditioned primal-dual algorithm with a nonlinear and asymmetric preconditioner to converge. 
Similar algorithms based on the separate and project principle for variational inequalities can refer to \cite{konnov2005combined,konnov2006combined,giselsson2021nonlinear,bui2020warped}.
We propose the nonlinear preconditioned primal-dual with projection method, which is described in Algorithm \ref{alg:1}.

\begin{algorithm}[H]
	\renewcommand{\algorithmicrequire}{\textbf{Input:}}
	\renewcommand{\algorithmicensure}{\textbf{Output:}}
	\caption{Nonlinear preconditioned primal-dual with projection}
	\label{alg:1}
	\begin{algorithmic}[1]
		\REQUIRE $u_0 \in\mathcal{H}$
		\ENSURE $u_N$
	    \FOR{$k=0,1,\cdots$}
		\STATE $r_{k}=(M+P)^{-1}Mu_{k}$
		\STATE $H_k=\{z:\langle z-r_k,(M+B)u_k-(M+B)r_k\rangle
		\leq \frac{L}{4}\|u_k-r_k\|^2\}$
		\STATE $u_{k+1}=(1-\theta_k)u_k+\theta_k\Pi_{H_k}(u_k)$
		\ENDFOR
	\end{algorithmic}  
\end{algorithm}

The first step can be viewed as a prediction step of (\ref{3:4}) with the preconditioner $M$ (\ref{3:1}).
The second step construct a halfspace $H_k$ between the current iterate $u_k$ and the solution set $\operatorname{zer}P$, which has been strictly verified. 
The iterate $u_{k+1}$ is in the third step obtained by a relaxed projection from $u_k$ onto the separating halfspace $H_k$. The relaxed parameter is limited as $\theta_k\in (0,2)$, and it determines the position of the update $u_{k+1}$ on an open segment between the current iterate $u_k$ and $2\Pi_{H_k}(u_k)-u_k$. We can accelerate the iterates toward a solution by selecting a proper relaxation sequence $\{\theta_k\}_{k\in\mathbb{\mathbb{N}}}$. Generally, we use a constant relaxation parameter $\theta_k= \theta$.

Here, we deduce the explicit form for the projection step in Algorithm \ref{alg:1}. Let $(M+B)u_k-(M+B)r_k=a_k$, representing the normal vector of the corresponding hyperplane. 
For all $u_k\notin H_k$, the projection of the current iterate onto the separating hyperspace is
$$\Pi_{H_k}(u_k)=u_k-\frac{\langle u_k-r_k,a_k\rangle- \frac{L}{4}\|u_k-r_k\|^2}{\|a_k\|^2} a_k.
$$
Let
$$
t_k=\frac{\langle u_k-r_k,a_k\rangle- \frac{L}{4}\|u_k-r_k\|^2}{\|a_k\|^2},
$$
then the fourth step can be reformulated as 
\begin{equation}
\label{eq3}
u_{k+1}=u_k-\theta_kt_k(M+B)(u_k-r_k).
\end{equation}
We have the following bounds on the stepsize $t_k$.
\begin{lemma}
	\label{lemma3:1}
	Suppose that Assumption \ref{assum:1} and \ref{assum:2} hold and $\mu>\gamma+\frac{L}{4}$. Then for all $u_k\notin \operatorname{zer} P$, the stepsize $t_k$ satisfies
	\begin{equation}
	\label{eq4}
	0<\frac{\mu-\gamma-\frac{L}{4}}{L+q}
	\leq
	t_k
	\leq
	\frac{\frac{3}{4}L+q}{\mu-\gamma}
	\end{equation}

\end{lemma} 
Lemma \ref{lemma3:1} can be directly deduced by Assumption \ref{assum:2}, and hence the proof is omitted here.

\begin{remark}
	\label{remark3:1}
	When $r_k=u_k$, we have $u_k\in \operatorname{fix} R=\operatorname{zer} P$, implying that $u_k$ is a solution of the operator inclusion problem $0\in Pu$. Thus, we can select $$\operatorname{residual}(k)=\|u_{k}-r_{k}\|$$
	as the stopping rule of Algorithm \ref{alg:1}. When the residual reaches our preset accuracy, we stop the algorithm. This stopping rule avoids discussing whether $\langle u_k-r_k,a_k\rangle$ is 0 as \cite{bui2020warped}.
\end{remark}

\subsection{The selection of operators in the preconditioner}
\label{section4.2}
Let us turn attention to the selection of operators in our proposed nonlinear and asymmetric preconditioner.
Based on Assumption \ref{assum:2}, we provide illustrations of the preconditioner in (\ref{3:1}).

\begin{example}[differentiable operators]
	\label{example1}
	Let $w_1$ and $w_2$ be strongly-convex, smooth and differentiable, and let $N_1=\nabla w_1,N_2=\nabla w_2$. The update in (\ref{3:4}) is
	\begin{equation}
	\label{eq3:1}
	\left\{
	\begin{array}{l}
	x_{k+1}=(\partial f+\nabla w_1)^{-1}(\nabla w_1 (x_{k})-\nabla_x\phi(x_k,y_k))\\
	y_{k+1}=(\partial g+\nabla w_2)^{-1}(\nabla w_2 (y_k)+\nabla_y\phi(x_k,y_k)+Q(x_{k}-x_{k+1}))
	\end{array},
	\right.
	\end{equation}
	the corresponding optimal problem is 
	$$
	\left\{
	\begin{array}{l}
	x_{k+1}=\arg\min_x\{f(x)+D_{w_1}(x,x_k)+\langle \nabla_{x}\phi(x_k,y_k),x\rangle\}\\
	y_{k+1}=\arg\min_y\{g(y)+D_{w_2}(y,y_k)-\langle \nabla_y\phi(x_k,y_k)+Q(x_k-x_{k+1}), y \rangle\}
	\end{array}.
	\right.
	$$
	For brevity, let $N_1=\frac{1}{\tau}I,N_2=\frac{1}{\sigma}I$, and corresponding iterate is
	\begin{equation}
	\label{3:9}
	\left\{
	\begin{array}{l}
	x_{k+1}=\text{prox}_{\tau f}(x_{k}-\tau\nabla_x\phi(x_k,y_k))\\
	y_{k+1}=\text{prox}_{\sigma g}(y_k+\sigma\nabla_y\phi(x_k,y_k)+\sigma Q(x_{k}-x_{k+1}))
	\end{array}.
	\right.
	\end{equation}
	Considering the convex-concave saddle-point problem (\ref{1:2}), let $Q=-(\theta+1)K$, then we can recover the classical primal-dual method in \cite{chambolle2011first} from (\ref{3:9}), i.e., (\ref{1:3}). 
	
	In particular, if $f(x)=\delta_{C_1}(x), g(y)=\delta_{C_2}(y)$, where $C_1$ and $C_2$ are convex sets in Hilbert space and $\delta(\cdot)$ is an indicator function, then (\ref{3:9}) can be rewritten as
	$$
	\left\{
	\begin{array}{l}
	x_{k+1}=\Pi_{C_1}(x_{k}-\tau\nabla_x\phi(x_k,y_k))\\
	y_{k+1}=\Pi_{C_2}(y_k+\sigma\nabla_y\phi(x_k,y_k)+\sigma Q(x_{k}-x_{k+1}))
	\end{array}.
	\right.
	$$
\end{example}

\begin{example}[Non-differentiable operators]
	Let $w_1,w_2$ be strongly convex, smooth, and non-differentiable and $f,g$ be differentiable, and then $N_1=\partial w_1,N_2=\partial w_2$.
	The prediction step in Algorithm \ref{alg:1} updates as
	$$\left\{
	\begin{array}{l}
	x_{k+1}=(\nabla f+ \partial w_1)^{-1}(\partial w_1 (x_{k})-\nabla_x\phi(x_k,y_k))\\
	y_{k+1}=(\nabla g+\partial w_2)^{-1}(\partial w_2 (y_k)+\nabla_y\phi(x_k,y_k)+Q(x_{k}-x_{k+1}))
	\end{array},
	\right.$$
	and the corresponding algorithm is
	$$
	\left\{
	\begin{array}{l}
	x_{k+1}=\arg\min_x\{f(x)+D_{w_1}^{s_k}(x,x_k)+\langle \nabla_{x}\phi(x_k,y_k),x\rangle\}\\
	y_{k+1}=\arg\min_y\{g(y)+D_{w_2}^{l_k}(y,y_k)-\langle \nabla_y\phi(x_k,y_k)+Q(x_k-x_{k+1}), y \rangle\}\\
	s_{k+1}\in s_k-\nabla f(x_{k+1})-\nabla_{x}\phi(x_k,y_k)\\
	l_{k+1}\in l_k-\nabla g(y_{k+1})+\nabla_y\phi(x_k,y_k)+Q(x_k-x_{k+1})
	\end{array}.
	\right.
	$$
	
\end{example}

\subsection{Remarks about parameters}
\label{section4.3}
For the convenience of convergence analysis in the following section, here we focus on the positive definiteness of some relevant operators.
It follows from the range of $t_k$ in (\ref{eq4}) and $\theta_k\in (0,2)$ that we can use constant forms, i.e., $\theta_k= \theta,t_k=t$.

Let $C=M+B$ and we have 
\begin{equation}
\label{eq4:2}
C=\left[
\begin{array}{cc}
N_1 & 0 \\
Q & N_2
\end{array}
\right].
\end{equation}
Further, decompose the operator $M$ into $M=DC$.
Let $H=\frac{1}{\theta t}D,N=\theta tC$, we have $M=DC=HN$.
To simplify the notation for analysis, we study the positive definiteness of two operators as follows:
\begin{equation}
\label{eq3:4}
D=I-BC^{-1}=I-B\left[
\begin{array}{cc}
N_1^{-1} & 0 \\
-N_2^{-1}QN_1^{-1} & N_2^{-1}
\end{array}
\right], 
G=M+M^T-N^THN.
\end{equation}

\begin{lemma}
	\label{lemma4.2}
	Suppose that Assumption \ref{assum:1} and Assumption \ref{assum:2} hold, if 
	\begin{equation}
	\label{eq55:1}
	\mu>\gamma+L,~2\mu>\theta t(L+q)q,
	\end{equation}
	then for the operators $D$ and $G$ in (\ref{eq3:4}), we obtain that $D$ and $G$ are positive definite operators.
\end{lemma}

Motivated by the paper \cite{he2012convergence}, we provide a suggestive relaxation parameter $\theta$. For the convenience of analysis, we consider the case that the operator $P$ is monotone instead of weakly monotone.
According to the optimality conditions of the prediction step of Algorithm \ref{alg:1}, we obtain that 
$$\left\{
\begin{array}{l}
f(x)-f(\tilde{x}_k)+(x-\tilde{x}_k)^T(N_1 (\tilde{x}_k-x_k)+\nabla_x\phi(u_k))\geq 0\\
g(y)-g(\tilde{y}_k)+(y-\tilde{y}_k)^T(N_2(\tilde{y}_k-y_k)-\nabla_y\phi(u_k)-Q(x_k-\tilde{x}_k))\geq 0
\end{array},
\right.$$
where $r_k=(\tilde{x}_k,\tilde{y}_k)$. Thus,
\begin{equation}
\label{4:6}
l(u)-l(r_k)+(u-r_k)^T[Br_k+M(r_k-u_k)]\geq 0.
\end{equation}
The interested readers are referred to the survey papers \cite{he2012convergence,he2022generalized} for more detail.
The correction step in Algorithm 1 becomes 
\begin{equation}
\label{4:7}
u_{k+1}=u_k-\theta t_kC(u_k-r_k). 
\end{equation}
To determine an appropriate value of $\theta$, we denote
$$
q_k(\theta):=\|u_k-u^*\|_D^2-\|u_{k+1}-u^*\|_D^2,
$$
which measures the progress in $k$-th iterate. It follows from the unknown of the true solution $u^*$ that we are unable to maximize $q_k(\theta)$ directly. 

Let $u=r_k$ in (\ref{eq1:4}), $u=u^*$ in (\ref{4:6}), and then we have
\begin{equation}
\label{eq55:2}
l(r_k)-l(u^*)+(r_k-u^*)^TBu^*\geq 0,
\end{equation}
\begin{equation}
\label{eq55:3}
l(u^*)-l(r_k)+(u^*-r_k)^T[Br_k+M(r_k-u_k)]\geq 0.
\end{equation}
Adding (\ref{eq55:2}) and (\ref{eq55:3}) and using the monotonicity of $B$ implies that
$$(u^*-r_k)^TM(r_k-u_k)\geq 0.$$
Therefore,
$$
(u_k-u^*)^TM(u_k-r_k)\geq (r_k-u_k)^TM(r_k-u_k).
$$
We have
\begin{equation}
\begin{aligned}
q_k(\theta)
&=\|u_k-u^*\|_D^2-\|u_k-\theta t_k C(u_k-r_k)-u^*\|_D^2\\
&=2\theta t_k(u_k-u^*)^TDC(u_k-r_k)-\theta^2t_k^2\|  C(u_k-r_k)\|_D^2\\
&\geq 2\theta t_k (r_k-u_k)^TDC(r_k-u_k)-\theta^2t_k^2\|C(u_k-r_k)\|_D^2.
\end{aligned}
\end{equation}
The lower bound of $q_k(\theta)$ is a quadratic function of $\theta$ and reaches maximum at 
\begin{eqnarray}
\label{4:10}
\theta_k^*=\frac{(r_k-u_k)^TM(r_k-u_k)}{t_k\|C(u_k-r_k)\|_D^2}.
\end{eqnarray}
Hence, when solving the convex-concave saddle-point problem, we can choose $\theta_k^*$ as above for the relaxation parameter of the correction step in Algorithm \ref{alg:1}. In specific, uponing observing $C^TDC=C^TM$, we do not need to compute the operator $D$ in (\ref{4:10}).

\section{Convergence}
In this section, we will analyze the convergence of Algorithm \ref{alg:1} under different assumptions.

\subsection{Weak convergence}
We first present a general result on weak convergence of sequences to the solution set $\operatorname{zer} P$.

\begin{theorem}
	\label{th2}
	Suppose that Assumption \ref{assum:1} and Assumption \ref{assum:2} with $\mu>\gamma+\frac{L}{4}$ hold and that the relaxation parameter $\theta_k\in (0,2)$ satisfies $\lim\inf_{k\rightarrow\infty} \theta_k(2-\theta_k)>0$. Then sequences $\{u_k\}_{k\in \mathbb{N}}$ and $\{r_k\}_{k\in \mathbb{N}}$ generated by Algorithm 1 such that
	\begin{enumerate}[(i)]
		\item $\{\|u_k-z\|\}_{k\in \mathbb{N}}$ converges for every $z\in \operatorname{zer} P$.
		\item $\{\|u_k-r_k\|\}_{k\in \mathbb{N}}$ converges to 0.
		\item $\{u_k\}_{k\in \mathbb{N}}$ converges weakly to a point in $\operatorname{zer} P$, i.e., $u_k\rightharpoonup \bar{u}\in \operatorname{zer} P$ as $k\rightarrow \infty$.
	\end{enumerate}

\end{theorem}

\begin{remark}
	Theorem \ref{th2} (ii) shows that $u_k$  persistently approximates $r_k$, which verifies that the stopping rule introduced in Remark \ref{remark3:1} is reasonable. 
\end{remark}

\subsection{Sublinear convergence}
\label{sec4:3}
We assume $P$ to be monotone, which can be realized by assuming that $\phi$ is convex-concave.
Recall notations in Section \ref{section4.3}, Algorithm \ref{alg:1} can be reformulated as
\begin{equation}
\label{eq4:1}
\begin{array}{c}
l(u)-l(r_k)+(u-r_k)^T[Br_k+M(r_k-u_k)]\geq 0,\\
u_{k+1}=u_k-\theta_k t_k C(u_k-r_k).
\end{array}
\end{equation}  
To simplify the convergence analysis, we use the constant parameter $\theta\in (0,2)$ and constant stepsize $t$ satisfying (\ref{eq4}).
Next, we derive the strict contraction of the sequence $\{u_k\}_{k\in \mathbb{N}}$ generated by Algorithm \ref{alg:1}, and then the $O(\frac{1}{n})$ convergence rate can be derived. 

\begin{theorem}
	\label{th4}
	Suppose that Assumption \ref{assum:1} and Assumption \ref{assum:2} hold with $\mu>\gamma+\frac{L}{4}$. In addition, assume that $P$ is monotone.
	The sequences $\{u_k\}_{k\in \mathbb{N}}$ and $\{r_k\}_{k\in \mathbb{N}}$ are generated by Algorithm \ref{alg:1}. Under the condition (\ref{eq55:1}) that $D$ and $G$ are positive definite, it holds that
	$$
	\|u_{k+1}-u^*\|_H^2\leq \|u_{k}-u^*\|_H^2-\|u_k-r_k\|_G^2.
	$$
\end{theorem}

Now, we are ready to prove the global convergence of the preconditioned primal-dual with projection in the case that $P$ is monotone.
\begin{theorem}
	\label{th5}
	Suppose that Assumption \ref{assum:1} and Assumption \ref{assum:2} hold with $\mu>\gamma+\frac{L}{4}$ and $P$ is monotone.
	Let $\{u_k\}_{k\in \mathbb{N}}$ and $\{r_k\}_{k\in \mathbb{N}}$ be the sequence generated by Algorithm \ref{alg:1}. Then, under the condition (\ref{eq55:1}) that $D$ and $G$ are positive definite, for any positive integer $n$, we have  $$\|N(u_{n+1}-r_{n+1})\|_H^2<\|N(u_{n}-r_{n})\|_H^2.$$ 
	and 
	\begin{equation}
	\label{eq4:3}
	\|N(u_n-r_n)\|_H^2\leq\frac{1}{(n+1)c} \|u_0-u^*\|_H^2,
	\end{equation}
	where $c>0$ is a constant independent of $n$.	
\end{theorem}
Recalling the reformulation (\ref{eq4:1}) of Algorithm \ref{alg:1}, 
we have $u_k-u_{k+1}=N(u_k-r_k)$;
we find that $r_k$ is a solution of the VI (\ref{eq1:4}) iff $\|N(u_k-r_k)\|^2_H=0$; the assertion (\ref{eq4:3}) implies that both $\|u_k-r_k\|<\varepsilon$ and $\|u_k-u_{k+1}\|<\varepsilon$ can be used as a stopping criterion, where $\varepsilon>0$ denotes the error tolerance. Theorem \ref{th5} also indicates the $O(\frac{1}{n})$ convergence rate for Algorithm \ref{alg:1}.

\subsection{Linear convergence}
\label{sec4:2}
Here, we show local linear convergence of Algorithm \ref{alg:1} under metric subregularity in \cite{dontchev2004regularity}. The distance between a point $x$ and a convex set $X$ is defined by
$$\operatorname{dist}(x,X):=\min\{\|x-y\|:y\in X\}.$$
The metric subregularity is defined as follows.
\begin{definition}[metric subregularity]
	A map $F:\mathcal{H}\rightarrow2^{\mathcal{H}}$ is metrically subregular at $\bar{x}$ for $\bar{y}$, if $\bar{y}\in F\bar{x}$ and there exists $\kappa\in[0,+\infty)$ along with neighborhoods $U$ of $\bar{x}$ and $V$ of $\bar{y}$ such that
	$$
	\operatorname{dist}(x,F^{-1}\bar{y})
	\leq \kappa \operatorname{dist}(\bar{y}, Fx\cap V),\forall x\in U.
	$$
	The infimum of the set of values $\kappa$ for which this holds is the modulus of metric
	subregularity, denoted by subreg $F(\bar{x}|\bar{y})$.
\end{definition}

Metric subregularity is weaker than metric regularity. We refer the reader to \cite{dontchev2004regularity,dontchev2009implicit} for more details on metric subregularity and its relation with the other regularity properties. A similar convergence result under metric subregularity is provided in \cite{latafat2017asymmetric}, and our proof is similar to Theorem 3.3 in \cite{latafat2017asymmetric}.
\begin{theorem}
	\label{th3}
	Suppose that Assumption \ref{assum:1} and Assumption \ref{assum:2} with $\mu>\gamma+\frac{L}{4}$ hold and the relaxation parameter $\theta_k\in[\varepsilon_{\theta},2-\varepsilon_{\theta}]$ for some $\varepsilon_{\theta}\in(0,1)$. Further assume that $P$ is metrically subregular at all $z\in \operatorname{zer}P$ for 0 with $\kappa \geq 0$, and $\mathcal{H}$ is  either finite-dimentional or the neighborhood $\mathcal{U}=\mathcal{H}$ at all $z\in \operatorname{zer}P$. Then, $\{u_k\}_{k\in \mathbb{N}}$ generated by Algorithm \ref{alg:1} converges to the solution set, i.e., $\operatorname{dist}(u_k,\operatorname{zer}P)$ converges locally Q-linearly to 0
	and $u_k\rightarrow \bar{u}\in \operatorname{zer} P$ locally R-linearly.
\end{theorem}

\begin{remark}
	Employing the constant parameter $\theta\in(0,2)$ and constant stepsize $t$ satisfying (\ref{eq4}), Theorem \ref{th2} and \ref{th3} still hold.
\end{remark}

\section{Connections to existing algorithms}
Here, we consider a saddle-point problem with a bilinear term (\ref{1:2}), i.e., $\phi(x,y)=\langle Kx,y\rangle$ with a continuous linear operator $K$.
We propose a relaxed preconditioned primal-dual with projection method and establish connections between our proposed framework and other well-known algorithms.
We find that our framework can recover many lastest primal-dual algorithms, e.g., \cite{zhu2008,chambolle2011first,he2022generalized}.
Instead of using the exact correction step, we propose a relaxed version of our framework as follows. 
\begin{algorithm}[H]
	\renewcommand{\algorithmicrequire}{\textbf{Input:}}
	\renewcommand{\algorithmicensure}{\textbf{Output:}}
	\caption{Relaxed preconditioned primal-dual with projection}
	\label{alg:2}
	\begin{algorithmic}[1]
		\REQUIRE $u_0 \in\mathcal{H}$
		\ENSURE $u_N$
		\FOR{$k=0,1,\cdots$}
		\STATE $r_{k}=Ru_{k}$
		\STATE $u_{k+1}=u_k-G(u_k-r_k)$
		\ENDFOR
	\end{algorithmic}  
\end{algorithm}

In Algorithm \ref{alg:2}, the relaxed operator $G$ requires to be strongly monotone and maixmally monotone; the corresponding iterates can be represented as the following steps.
\begin{equation}
\label{eq6:3}
\left\{
\begin{array}{l}
\tilde{x}_k=(N_1+\partial f)^{-1}(N_1x_{k}- K^Ty_k)\\
\tilde{y}_k=(N_2+\partial g)^{-1}(N_2y_k-Q\tilde{x}_{k}+(K+Q)x_k)\\
u_{k+1}=u_k-G(u_k-r_k)
\end{array},
\right.
\end{equation}
where where $r_k=(\tilde{x}_k,\tilde{y}_k)$.

Our relaxed framework can recover the classical primal-dual method \cite{chambolle2011first} by letting $$N_1=\frac{1}{\tau}I,N_2=\frac{1}{\sigma}I,Q=-(\theta+1)K,G=I.$$
Especially, let $\theta=1$, our algorithm reduces to PDHG method \cite{zhu2008}, updating as
\begin{equation}
\left\{
\begin{array}{l}
x_{k+1}=\text{prox}_{\tau f}(x_k- \tau K^Ty_k)\\
y_{k+1}=\text{prox}_{\sigma g}(y_k+\sigma K(2x_{k+1}-x_k))
\end{array}.
\right.
\end{equation}
The generalized primal-dual method proposed by \cite{he2022generalized} can be written as the precondition-correction representation
\begin{equation}
\label{eq6:2}
\left\{\begin{array}{l}
\tilde{x}_k=\text{prox}_{\tau f}( x_k- \tau K^Ty_k)\\
\tilde{y}_k=\text{prox}_{\sigma g}( y_k+\sigma(1+\theta)K\tilde{x}_k-\sigma\theta Kx_k)))\\
u_{k+1}=u_k-W(u_k-r_k)
\end{array},
\right.
\end{equation}
where $\theta\in [0,1]$ and
$$W=\left[
\begin{array}{cc}
I & 0 \\
(1-\theta)\sigma K & I
\end{array}
\right].$$
Algorithm \ref{alg:2} can recover generalized primal-dual algorithm in \cite{he2022generalized} by letting 
\begin{equation}
\label{eq6:4}
N_1=\frac{1}{\tau}I,N_2=\frac{1}{\sigma}I,Q=-(\theta+1)K,G=W.
\end{equation}
Recall the main convergence result in Theorem 5.3 in \cite{he2022generalized}
$$
\|G(u_n-r_n)\|_H^2\leq \frac{1}{(n+1)c}\|u_0-u^*\|_H^2.
$$
In the case of (\ref{eq6:4}), we find that
$$
G=\left[
\begin{array}{cc}
I & 0 \\
(1-\theta)\sigma K & I
\end{array}
\right],
M=
\left[
\begin{array}{cc}
\frac{1}{\tau}I & -K^T \\
-\theta K & \frac{1}{\sigma}I
\end{array}
\right],
H=MG^{-1}=
\left[
\begin{array}{cc}
\frac{1}{\tau} I+(1-\theta)\sigma K^T K & -K^T \\
-K & \frac{1}{\sigma} I
\end{array}
\right].
$$
Therefore, Theorem \ref{th5} in our paper recovers the convergence result of the generalized primal-dual algorithm. 

\begin{remark}
For bilinear saddle-point problems (\ref{1:2}), primal-dual method with preconditioner $M$ in (\ref{3:1}) can recover existing algorithms, e.g., \cite{zhu2008,chambolle2011first,he2022generalized}; For nonlinear saddle-point problems (\ref{eq1:1}), preconditioner $M$ in (\ref{3:2}) provides a new and general form. 
\end{remark}




\section{Conclusions}
In this paper, we present the general algorithm nonlinear preconditioned primal-dual with projection for solving nonconvex-nonconcave and non-smooth saddle-point problems. It is based on a nonlinear and asymmetric map called warped resolvent and belongs to the separate and project method. The warped mapping step constructs a separating hyperspace between the current point and the solution set. The next step is to project onto the hyperspace to further accelerate the original warped proximal iterations and ensure the convergence of the warped preconditioned step with nonlinear and asymmetric preconditioned operators. We prove weak convergence to the solution set, sublinear convergence under the convexity of the saddle-point problems, and linear convergence under a metric subregularity assumption. We also prove the generalized primal-dual algorithm and NOFOB method are special cases of our algorithm. 

In the future, we would like to consider the form of preconditioner $M$ in (\ref{3:2}) and explore the convergence properties of (\ref{3:4}) with correction (\ref{eq2:14}).
Meanwhile, the advantages and disadvantages of the two forms of preconditioner are also worth analyzing.

\section*{Acknowledgements}
This work was supported by the National Science Foundation of China (No.11971480),
the Natural Science Fund of Hunan for Excellent Youth (No.2020JJ3038), and the Fund for
NUDT Young Innovator Awards (No.20190105).

\section*{Appendices}
\begin{appendices}
	\textbf{Proof of lemma \ref{lemma:2}.} 
	\begin{enumerate}[(i)]
		\item Let $u=(x,y),y=(\bar{x},\bar{y})$, note that $\|\nabla \phi (x,y)-\nabla \phi (\bar{x},\bar{y})\|=\|Bu-Bv\|$. Since
		$$
		\begin{aligned}
		\|\nabla \phi (x,y)-\nabla \phi (\bar{x},\bar{y})\|
		&\leq
		\|\nabla \phi (x,y)-\nabla \phi (\bar{x},y)\|
		+
		\|\nabla \phi (\bar{x},y)-\nabla \phi (\bar{x},\bar{y})\| \\
		&\leq \sqrt{L_{xx}^2+L_{yx}^2}\|x-\bar{x}\|+\sqrt{L_{xy}^2+L_{yy}^2}\|y-\bar{y}\|\\
		&\leq \sqrt{L_{xx}^2+L_{yx}^2+L_{xy}^2+L_{yy}^2}\|u-v\|\\
		&=L\|u-v\|,
		\end{aligned}
		$$  
		namely, $\|Bu-Bv\|\leq L\|u-v\|$, we say $B$ is $L$-Lipschitz continuous. Recall that
		$\phi(\cdot, y)$ is $L_{x x}$-weakly convex for fixed $y$, and $-\phi(x, \cdot)$ is $L_{y}$-weakly convex for fixed $x$. Then, using the weak convexity on $x$, we have
		$$
		\begin{aligned}
		& \phi(\overline{x}, y) \geq \phi(x, y)+\left\langle\nabla_{x} \phi(x, y), \overline{x}-x\right\rangle-\frac{L_{x x}}{2}\|\overline{x}-x\|^2, \\
		& \phi(x, \overline{y}) \geq \phi(\overline{x}, \overline{y})+\left\langle\nabla_{x} \phi(\overline{x}, \overline{y}), x-\overline{x}\right\rangle-\frac{L_{x x}}{2}\|x-\overline{x}\|^2,
		\end{aligned}
		$$
		for all $x, \overline{x} \in \mathcal{H}$ and $y, \overline{y} \in \mathcal{H}$. Similarly, using the weak convexity on $y$, we have
		$$
		\begin{aligned}
		& -\phi(x, \overline{y}) \geq-\phi(x, y)-\left\langle\nabla_{y} \phi(x, y), \overline{y}-y\right\rangle-\frac{L_{y y}}{2}\|\overline{y}-y\|^2, \\
		& -\phi(\overline{x}, y) \geq-\phi(\overline{x}, \overline{y})-\left\langle\nabla_{y} \phi(\overline{x}, \overline{y}), y-\overline{y}\right\rangle-\frac{L_{y y}}{2}\|y-\overline{y}\|^2 .
		\end{aligned}
		$$
		Then, summing the above four inequalities yields
		$$
		\begin{aligned}
		~~~~ \langle u-v, Bu-B v\rangle &\geq-L_{x x}\|x-\overline{x}\|^2-L_{y y}\|y-\overline{y}\|^2 \\
		& \geq-\max \left\{L_{x x}, L_{y y}\right\}\|x-y\|^2
		=
		-\gamma \|u-v\|^2,
		\end{aligned}
		$$
		thus, $B$ is $\gamma$-weakly monotone.
		\item Recall that $P=A+B$, then for all $u=(x,y),v=(\bar{x},\bar{y})\in \mathcal{H}$,
		$$
		\begin{aligned}
		&~~~~\langle Pu-Pv,u-v\rangle\\
		&=(x-\bar{x})^T(\partial f(x)-\partial f(\bar{x}))+(y-\bar{y})^T(\partial g(y)-\partial g(\bar{y}))+\langle B u-Bv,u-v\rangle\\
		& \geq \langle Bu-Bv,u-v\rangle\\
		&\geq -\gamma \|u-v\|^2.
		\end{aligned}
		$$
		The first inequality holds with the convexity of $f$ and $g$. The second inequality uses the $\gamma$-weakly monotonicity of $B$.
	\end{enumerate}
	\textbf{Proof of Lemma \ref{lemma:3}.} 
For $u=(x,y),v=(\bar{x},\bar{y})$,
	$$
	\begin{aligned}
	&~~~~\langle Mu-Mv,u-v\rangle\\
	&=
	-\langle Bu-Bv,u-v\rangle
	+\langle x-\bar{x},N_1 x-N_1 \bar{x}\rangle +\langle x-\bar{x},Q y-Q\bar{y}\rangle+\langle y-\bar{y},N_2y-N_2\bar{y}\rangle\\
	&\geq -L\|u-v\|^2+\mu_1 \|x-\bar{x}\|^2+\mu_2 \|y-\bar{y}\|^2-\frac{\|Q\|}{2}\|x-\bar{x}\|^2-\frac{\|Q\|}{2} \| y- \bar{y}\|^2\\
	&\geq(\mu_3-L-\frac{\|Q\|}{2})\|u-v\|^2,\\
	&=\mu\|u-v\|^2,
	\end{aligned}
	$$
	where $\mu=\mu_3-L-\frac{\|Q\|}{2},\mu_3=\min\{\mu_1,\mu_2\}$. If $\mu>0$, then $M$ is $\mu$-strongly monotone. From
	$$
	\begin{aligned}
	&~~~~\langle Mu-Mv,u-v\rangle\\
	&=
	-\langle Bu-Bv,u-v\rangle
	+\langle x-\bar{x},N_1 x-N_1 \bar{x}\rangle +\langle x-\bar{x},Q y-Q\bar{y}\rangle+\langle y-\bar{y},N_2y-N_2\bar{y}\rangle\\
	&\leq \gamma\|u-v\|^2+L_1 \|x-\bar{x}\|^2+L_2\|y-\bar{y}\|^2+\frac{\|Q\|}{2}\|x-\bar{x}\|^2+\frac{\|Q\|}{2} \| y- \bar{y}\|^2\\
	&\leq(\gamma+L'+\frac{\|Q\|}{2})\|u-v\|^2,\\
	&=q\|u-v\|^2,
	\end{aligned}
	$$
	where $q=\gamma+L'+\frac{\|Q\|}{2},L'=\max\{L_1,L_2\}$, we conclude that $M$ is $q$-Lipschitz continuous. 
	\\
	\textbf{Proof of Lemma \ref{lemma:4}.} 
	\begin{enumerate}[(i)]
		\item 
		Note that $P+M=(P+\gamma I)+(M-\gamma I)$. From the condition that $\mu>\gamma$ and $P$ is $\gamma$-weakly maximally monotone, it is straightforward to show that $P+\gamma I$ is maximally monotone and $M-\gamma I$ is strongly monotone. Then \cite{bauschke2008general} shows that $\operatorname{ran}(P+M)=\mathcal{H}$. This implies that $Rz$ is nonempty for all $z\in\mathcal{H}$. Assume that $u,v\in Rz$. Since $Pu\in Mz-Mu,Pv\in Mz-Mv$, we have
		$$
		\begin{aligned}
		-\gamma\|u-v\|^2
		&\leq 
		-\langle Mu-Mv,u-v\rangle\\
		&=
		\langle Mz-Mu,u-v\rangle-\langle Mz-Mv,u-v\rangle\\
		&=\langle Pu-Pv,u-v\rangle\\
		&\leq 
		-\mu\|u-v\|^2.
		\end{aligned}
		$$
		If $\mu >\gamma$, the above inequality implies that $u=v$.
		\item Applying $\operatorname{dom}(M+P)^{-1}=\operatorname{ran}(M+P)=\mathcal{H}$, hence $R$ has full domain.
		\item $u\in \operatorname{fix}R$ is equivalent to
		$$
		u=(M+P)^{-1}Mu\Leftrightarrow Mu\in Mu+Pu\Leftrightarrow 0\in Pu.
		$$
		This concludes the proof. 
	\end{enumerate}
\textbf{Proof of Theorem \ref{th1}.} 
\begin{enumerate}[(i)]
		\item According to the $\mu$-strong monotonicity of $M$, the $\gamma$-weakly monotone of $B$, and $\mu>\gamma$, we have
		$$\begin{aligned}
		\psi_u(u)&=\langle (M+B)u-(M+B)Ru,u-Ru\rangle-\frac{L}{4}\|u-Ru\|^2.\\
		&=\langle Mu-MRu,u-Ru\rangle+\langle Bu-BRu,u-Ru\rangle-\frac{L}{4}\|u-Ru\|^2.\\
		&\geq (\mu-\gamma-\frac{L}{4})\|u-Ru\|^2.
		\end{aligned}$$
		This implies that $\psi_u(u)\geq 0$ for all $u\in\mathcal{H}$ when $\mu>\gamma+\frac{L}{4}$.
		\item That $\psi_u(u)>0$ for all $u\notin \operatorname{fix}R$ from (i). If $u\in \operatorname{fix}R$, then $Ru=u,\psi_u(u)=0$.
		\item
	    For all $z\in \operatorname{zer}P, u\in\mathcal{H}$, in subsection \ref{section4:1} we have proved that 
		$$
		\begin{aligned}
		\langle (M+B)u-(M+B)Ru,z-Ru\rangle-\frac{L}{4}\|u-Ru\|^2,\leq 0
		\end{aligned}
		$$
	Thus, $\psi_u(z)\leq 0$.
	\end{enumerate}
\textbf{Proof of Lemma \ref{lemma4.2}.}
	For all $u,v\in\mathcal{H}$, we have
	$$
	\langle u,Dv\rangle=\langle Du,v\rangle,
	$$
	and
	\begin{equation*}
	\begin{aligned}
	(u-v)^TD(u-v)
	&=\|u-v\|^2-(u-v)^TBC^{-1}(u-v)\\
	&\geq \|u-v\|^2-\|B(u-v)\|\cdot\|C^{-1}(u-v)\|\\
	&\geq (1-\frac{L}{\mu-\gamma})\|u-v\|^2.
	\end{aligned}
	\end{equation*}
	The first inequality uses Cauchy-schwartz inequality, and the second inequality holds with the $L$-Lipschitz continuity of $B$ and $(\mu-\gamma)$-strong monotonicity of $C$. 
	Thus, the operator $D$ is a positive definite operator when $\mu>\gamma+L$. 
	
	It is clear that 
	\begin{equation*}
	\langle u,Gv\rangle = \langle Gu,v\rangle,\forall u,v\in\mathcal{H}
	\end{equation*}
	and
	\begin{equation*}
	\begin{aligned}
	(u-v)^TG(u-v)
	&=2(u-v)^T(Mu-Mv)-(u-v)^TN^TM(u-v)\\
	&\geq 2\mu\|u-v\|^2-\theta t\|C(u-v)\|\cdot\|M(u-v)\|\\
	&\geq (2\mu-\theta t(L+q)q)\|u-v\|^2.
	\end{aligned}
	\end{equation*}\\
	The first inequality holds with the $\mu$-strong monotonicity of $M$ and Cauchy-schwartz inequality. The second inequality uses the Lipschitz continuity of $C$ and $M$. 
	Therefore, the operator $G$ is positive definite when $2\mu>\theta t(L+q)q$.\\
\textbf{Proof of Theorem \ref{th2}.} 
\begin{enumerate}[(i)]
	\item 
	For all $z\in \operatorname{zer} P$, according to the projection theorem, we have
	\begin{equation}
	\label{eq5:1}
	\begin{aligned}
	\|u_{k+1}-z\|^2
	&=\|u_k-z+\theta_k(\Pi_{H_k}u_k-u_k)\|^2\\
	&=\|u_k-z\|^2+\theta_k^2\|\Pi_{H_k}u_k-u_k\|^2+2\theta_k\langle u_k-z,\Pi_{H_k}u_k-u_k\rangle\\
	&\leq \|u_k-z\|^2-\theta_k(2-\theta_k)\|\Pi_{H_k}u_k-u_k\|^2,
	\end{aligned}
	\end{equation}
	where $\theta_k\in (0,2)$.
Thus, $\{u_k\}_{k\in\mathbb{N}}$ is bounded and $\{\|u_k-z\|\}_{k\in \mathbb{N}}$ is nonincreasing, lower bound, i.e., it converges.
	\item Summing (\ref{eq5:1}), it implies 
	$$\sum_{k\in\mathbb{N}}\theta_k(2-\theta_k)\|\Pi_{H_k}u_k-u_k\|^2\leq \|u_0-z\|^2<+\infty.$$
	Since the parameter $\theta_k\in (0,2)$ satisfies $\lim\inf_{k\rightarrow\infty} \theta_k(2-\theta_k)>0$, we have $\|\Pi_{H_k}u_k-u_k\|\rightarrow 0$ as $k\rightarrow \infty$.
	From Theorem \ref{th1}, since $u_k\notin \operatorname{zer}P$, we have $u_k\notin H_k$. By the $(\mu_3-\frac{\|Q\|}{2})$-strongly monotonicity of $M+B$ and $t_k\geq \frac{\mu-\gamma-\frac{L}{4}}{L+q}>0$, we obtain
	\begin{equation}
	\label{eq5:2}
	0\leftarrow \|\Pi_{H_k}u_k-u_k\|=t_k\|(M+B)u_k-(M+B)r_k\|\geq \frac{2\mu_3-\|Q\|}{2}\cdot\frac{\mu-\gamma-\frac{L}{4}}{L+q}\|u_k-r_k\|\geq 0,
	\end{equation}
	and hence $\|u_k-r_k\|\rightarrow 0$.
	\item  By Lemma 2.47 in \cite{bauschke2011convex} and (i), it is enough to show that every weak sequential cluster point belongs to $\operatorname{zer}P$. By Lemma 2.45 in \cite{bauschke2011convex}, at least one cluster point exists since (i) implies that $\{u_k\}_{k\in \mathbb{N}}$ is bounded. \\
	Let $\{u_{n_k}\}_{k\in \mathbb{N}}$ be a weakly convergent subsequence $u_{n_k}\rightharpoonup u$, where $u$ is the cluster point. Since $r_{n_k}=Ru_{n_k}=(M+P)^{-1}Mu_{n_k}$, it holds that
	$$v_{n_k}:=Mu_{n_k}-Mr_{n_k}\in Pr_{n_k}.$$
	Since $u_{n_k}\rightharpoonup u$ and $u_{n_k}-r_{n_k}\rightarrow 0$, we conclude $r_{n_k}\rightharpoonup u$.
	Further, by $q$-Lipschitz continuity of $M$, we conclude that
	$$\|v_{n_k}\|=\|Mu_{n_k}-Mr_{n_k}\|\leq q\|u_{n_k}-r_{n_k}\|\rightarrow0,$$
	thus, $v_{n_k}\rightarrow 0$. 
	Since $P$ is maximally monotone and $(r_{n_k},v_{n_k})\in \operatorname{gra} P$, we conclude by Proposition 20.38 in \cite{bauschke2011convex} that $(u,0)\in \operatorname{gra} P$, i.e., $u\in \operatorname{zer}P$. To sum up, every weak sequential cluster point belongs to $\operatorname{zer} P$. The proof is completed.
\end{enumerate}
\textbf{Proof of Theorem \ref{th3}.}
	Observing $P^{-1}(0)=\operatorname{zer} P$,  $P$ is metrically subregular at all $z\in \operatorname{zer}P$ for 0 implies that there exists $\kappa\geq 0$ along with neighborhoods $\mathcal{U}$ of $u$ and $\mathcal{V}$ of 0 such that
	\begin{equation}
	\label{eq5:3}
	\operatorname{dist}(u,\operatorname{zer} P)\leq \kappa \operatorname{dist}(0,Pu\cap \mathcal{V})\leq \kappa \|v\|
	\end{equation}
	for all $u\in \mathcal{U}, v\in \mathcal{V}, v\in Pu$. Recall that $v_k=Mu_k-Mr_k\in Pr_k$ and $r_k=Ru_k$.
	By the $q$-Lipschitz continuty of $M$ and Theorem \ref{th2}, we have
	$$
	\|v_k\|=\|Mu_k-Mr_k\|\leq q\|u_k-r_k\|\rightarrow 0.
	$$
	Therefore, there exists a $K_v\in N$ such that $v_k\in Pr_k\cap \mathcal{V}$. 
	If $\mathcal{U}=\mathcal{H
	}$, we have $r_k\in\mathcal{U}$ for all $k\in \mathbb{N}$. If $\mathcal{H}$ is finite-dimensional, by $u_k-r_k\rightarrow0$ and $u_k\rightharpoonup \bar{u}\in \operatorname{zer} P$ in Theorem 
	\ref{th2} and a neighborhood of $\operatorname{zer} P$ is $\mathcal{U}$, we have there exists a $K_{r}\in N$ such that $r_k\in\mathcal{U}$ for all $k\geq K_r$.
	Therefore, from (\ref{eq5:3}) for all $k\geq K:=\max\{K_v,K_r\}$ we have
	$$
	\operatorname{dist}(r_k,\operatorname{zer}P)\leq\kappa \operatorname{dist}(0,Pr_k\cap \mathcal{V})\leq    \kappa \|v_k\|\leq \kappa q\|u_k-r_k\|,
	$$
	thus,
	\begin{equation}
	\label{eq5:4}
	\begin{aligned}
	\operatorname{dist}(u_k,\operatorname{zer}P)
	&=
	\|u_k-\Pi_{\operatorname{zer} P}u_k\|\\
	&\leq \|u_k-\Pi_{\operatorname{zer} P}r_k\|\\
	&\leq \|u_k-r_k\|+\|r_k-\Pi_{\operatorname{zer} P}r_k\|\\
	&=\|u_k-r_k\|+\operatorname{dist}(r_k,\operatorname{zer} P)\\
	&\leq (1+\kappa q)\|u_k-r_k\|
	\end{aligned}
	\end{equation}
	for all $k\geq K:=\max\{K_v,K_r\}$. 
	Let $z=\Pi_{\operatorname{zer} P}u_k,\theta_k\in[\varepsilon_{\theta},2-\varepsilon_{\theta}]$ for some $\varepsilon_{\theta}\in (0,1)$, and use (\ref{eq5:1}), (\ref{eq5:4}), and that followed by (\ref{eq5:2}) with $\mathcal{Z}:=\frac{2\mu_3-\|Q\|}{2}\cdot\frac{\mu-\gamma-\frac{L}{4}}{L+q}>0$ to obtain that
	\begin{equation}
	\label{eq5:5}
	\begin{aligned}
	\operatorname{dist}^2(u_{k+1},\operatorname{zer} P)&\leq
	\|u_{k+1}-z\|^2\\
	&\leq \|u_k-z\|^2-\theta_k(2-\theta_k)\|u_k-\Pi_{H_k}u_k\|^2\\
	&\leq \|u_k-z\|^2-\varepsilon_{\theta}(2-\varepsilon_{\theta})\|u_k-\Pi_{H_k}u_k\|^2\\
	&\leq \|u_k-z\|^2-\varepsilon_{\theta}(2-\varepsilon_{\theta})\mathcal{Z}^2\|u_k-r_k\|^2\\
	&\leq \left(1-\frac{\varepsilon_{\theta}(2-\varepsilon_{\theta})\mathcal{Z}^2}{(1+\kappa q)^2}\right)\operatorname{dist}^2(u_k,\operatorname{zer}P),
	\end{aligned}
	\end{equation}
	for all $k\geq K$.
	We obtain that $\operatorname{dist}(u_{k},\operatorname{zer} P)$ converges locally Q-linearly to 0 since $\varepsilon_{\theta}\in (0,1),\mathcal{Z}>0,\kappa\geq 0$ and $q>0$.
	
	Let 
	$$
	a:=\sqrt{1-\frac{\varepsilon_{\theta}(2-\varepsilon_{\theta})\mathcal{Z}^2}{(1+\kappa q)^2}},~d:=a^{-K}\sqrt{\frac{2-\varepsilon_{\theta}}{\varepsilon_{\theta}}}\operatorname{dist}(u_K,\operatorname{zer}P),
	$$
	where $a\in (0,1),R\geq 0$. 
	The updated step in Algorithm 1 implies that for all $k>K$
	$$
	\begin{aligned}
	\|u_{k+1}-u_k\|
	&=\theta_k\|u_k-\Pi_{H_k}u_k\|\\
	&\leq \sqrt{\frac{\theta_k}{2-\theta_k}}\operatorname{dist}(u_k,\operatorname{zer}P)\\
	&\leq \sqrt{\frac{2-\varepsilon_{\theta}}{\varepsilon_{\theta}}}\operatorname{dist}(u_k,\operatorname{zer}P)\\
	&\leq \sqrt{\frac{2-\varepsilon_{\theta}}{\varepsilon_{\theta}}}a^{k-K}\operatorname{dist}(u_K,\operatorname{zer}P)\\
	&=da^k,
	\end{aligned}
	$$
	the first inequality deduced by (\ref{eq5:1}), the second inequality holds with $\theta_k\in [\varepsilon_{\theta},2-\varepsilon_{\theta}]$, and the third inequality uses (\ref{eq5:5}). For all $j>k\geq K$, we have
	\begin{equation}
	\label{eq5:6}
	\|u_k-u_j\|\leq \sum_{i=0}^{j-1}\|u_{k+i+1}-u_{k+i}\|\leq d\sum_{i=0}^{j-1}c^{k+i}\leq dc^k\sum_{i=0}^{j-1}a^i=\frac{d}{1-a}a^k,
	\end{equation}
	so $\{u_k\}_{k\in \mathbb{N}}$ is a Cauchy sequence. Then $\{u_k\}_{k\in \mathbb{N}}$ converges. Let $u^*$ be the limit point, i.e., $u_k\rightarrow u^*,k\rightarrow \infty$. By (\ref{eq5:6}) we have
	\begin{equation}
	\label{eq5:7}
	\|u_k-u^*\|=\lim_{j\rightarrow\infty}\|u_k-u_j\|\leq \frac{d}{1-a}a^k.
	\end{equation}
	Therefore $\|u_k-u^*\|\rightarrow 0,k\rightarrow+\infty$. Moreover, by Theorem \ref{th2} we have
	$u_k\rightharpoonup \bar{u}\in \operatorname{zer}P$. We have $u^*=\bar{u},u_k\rightarrow \bar{u}\in \operatorname{zer} P$ from Corollary 2.52 in \cite{bauschke2011convex}. Given that $\{\frac{d}{1-a}a^k\}_{k\in \mathbb{N}}$ is linearly convergence and $\{u_k\}_{k\in \mathbb{N}}$ can be controlled by $\{\frac{d}{1-a}a^k\}_{k\in \mathbb{N}}$, we have $\{u_k\}_{k\in \mathbb{N}}$ is $R$-linearly convergence and the linear rate is given by (\ref{eq5:7}).\\
\textbf{Proof of Theorem \ref{th4}.}
According to (\ref{eq4:1}) we have
\begin{equation}
\label{eq5:8}
l(u)-l(r_k)+(u-r_k)^TBr_k\geq (u-r_k)^TH(u_k-u_{k+1}).
\end{equation} 
Applying the identity
$$
(a-b)^TH(c-d)=\frac{1}{2}\{\|a-d\|_H^2-\|a-c\|_H^2+\|c-b\|_H^2-\|d-b\|_H^2\}
$$
to the right-hand side of (\ref{eq5:8}), we have
\begin{equation}
\label{eq5:9}
(u-r_k)^TH(u_k-u_{k+1})
=\frac{1}{2}\{
\|u-u_{k+1}\|_H^2-\|u-u_k\|_H^2\}+\frac{1}{2}\{\|u_k-r_k\|_H^2-\|u_{k+1}-r_k\|_H^2
\}.
\end{equation}
For the second part of the right-hand side of (\ref{eq5:9}), because of $M=HN,2u^TMu=u^T(M^T+M)u$, we have
\begin{equation}
\label{eq7:1}
\begin{aligned}
\|u_k-r_k\|_H^2-\|u_{k+1}-r_k\|_H^2
&=
\|u_k-r_k\|_H^2-\|u_k-r_k-N(u_k-r_k)\|_H^2\\
&=
(u_k-r_k)^T(M+M^T-N^THN)(u_k-r_k)\\
&=\|u_k-r_k\|_G^2.
\end{aligned}
\end{equation}
Thus we prove the assertion from (\ref{eq5:8}), (\ref{eq5:9}) and (\ref{eq7:1})
\begin{equation}
\label{eq5:10}
l(u)-l(r_k)+(u-r_k)^TBr_k\geq
\frac{1}{2}\{
\|u-u_{k+1}\|_H^2-\|u-u_k\|_H^2+\|u_k-r_k\|_G^2
\}.
\end{equation}
Setting $u$ in (\ref{eq5:10}) as a solution $u^*\in \operatorname{zer}P$, we have
\begin{equation}
\label{eq7:2}
l(u^*)-l(r_k)+(u^*-r_k)^TBr_k\geq
\frac{1}{2}\{
\|u^*-u_{k+1}\|_H^2-\|u^*-u_k\|_H^2+\|u_k-r_k\|_G^2
\}.
\end{equation}
Let $u=r_k$ in (\ref{eq1:4}), implying
\begin{equation}
\label{eq7:3}
l(r_k)-l(u^*)+(r_k-u^*)^TBu^*\geq 0.
\end{equation}
Adding (\ref{eq7:2}) and (\ref{eq7:3}) and using the monotonicity of $B$, we deduce
$$
\|u_{k+1}-u^*\|_H^2\leq \|u_{k}-u^*\|_H^2-\|u_k-r_k\|_G^2.
$$
\textbf{Proof of Theorem \ref{th5}.} Using
$$
\|a\|_H^2-\|b\|_H^2=2a^TH(a-b)-\|a-b\|_H^2 
$$
with $a=N(u_k-r_k)$ and $b=N(u_{k+1}-r_{k+1})$, we obtain
\begin{equation}
\begin{aligned}
\label{eq5:11}
&~~~~\|N(u_k-r_k)\|_H^2-\|N(u_{k+1}-r_{k+1})\|_H^2\\
&=2(u_k-r_k)^TN^THN((u_k-r_k)-(u_{k+1}-r_{k+1}))-\|N((u_k-r_k)-(u_{k+1}-r_{k+1}))\|_H^2
\end{aligned}
\end{equation}
Firstly, we bound the first term of the right side in (\ref{eq5:11}). Setting $u=r_{k+1}$ in (\ref{eq4:1}), we have
\begin{equation}
\label{eq5:12}
l(r_{k+1})-l(r_k)+(r_{k+1}-r_k)^TBr_{k}\geq (r_k-r_{k+1})^TM(r_k-u_k).
\end{equation}
It follows from (\ref{eq4:1}) in $(k+1)$-th iteration that
\begin{equation}
\label{eq5:13}
l(u)-l(r_{k+1})+(u-r_{k+1})^TBr_{k+1}\geq (u-r_{k+1})^TM(u_{k+1}-r_{k+1}),\forall u\in \Omega.
\end{equation}
Setting $u=r_{k}$ in (\ref{eq5:13}), we get
\begin{equation}
\label{eq5:14}
l(r_{k})-l(r_{k+1})+(r_{k}-r_{k+1})^TBr_{k+1}\geq (r_{k}-r_{k+1})^TM(u_{k+1}-r_{k+1}).
\end{equation}
Adding (\ref{eq5:12}) and (\ref{eq5:14}), and using the monotonicity of $B$ we have
\begin{equation*}
\label{eq5:15}
0\geq (r_{k}-r_{k+1})^T(Br_{k+1}-Br_{k})\geq (r_{k}-r_{k+1})^TM\{-(r_{k+1}-u_{k+1})+(r_k-u_k)\},
\end{equation*}
thus,
\begin{equation}
\label{eq5:16}
(r_{k}-r_{k+1})^TM\{(r_{k+1}-u_{k+1})-(r_k-u_k)\}\geq 0.
\end{equation}
Adding $\{(r_{k+1}-u_{k+1})-(r_k-u_k)\}M\{(r_{k+1}-u_{k+1})-(r_k-u_k)\}$ to (\ref{eq5:16}), we have
\begin{equation}
\label{eq5:17}
(u_{k}-u_{k+1})^TM\{(r_{k+1}-u_{k+1})-(r_k-u_k)\}\geq 
\frac{1}{2}\|(r_{k+1}-u_{k+1})-(r_k-u_k)\|^2_{M+M^T}.
\end{equation}
Note that
\begin{equation*}   \begin{aligned}
&~~~~(u_{k}-u_{k+1})^TM\{(r_{k+1}-u_{k+1})-(r_k-u_k)\}\\
&=\{N(u_k-r_k)\}^TM\{(r_{k+1}-u_{k+1})-(r_k-u_k)\}\\
&=(u_k-r_k)^TN^THN\{(r_{k+1}-u_{k+1})-(r_k-u_k)\}.
\end{aligned}
\end{equation*}
Thus, (\ref{eq5:17}) can be rewritten as 
\begin{equation}
\label{eq5:18}
(u_k-r_k)^TN^THN\{(r_{k+1}-u_{k+1})-(r_k-u_k)\}\geq 
\frac{1}{2}\|(r_{k+1}-u_{k+1})-(r_k-u_k)\|^2_{M+M^T}.
\end{equation}
Substituting (\ref{eq5:18}) into (\ref{eq5:11}), we have
\begin{equation*}
\begin{aligned}
&~~~~\|N(u_k-r_k)\|_H^2-\|N(u_{k+1}-r_{k+1})\|_H^2\\
&\geq \|(r_{k+1}-u_{k+1})-(r_k-u_k)\|^2_{M+M^T}-\|N((u_k-r_k)-(u_{k+1}-r_{k+1}))\|_H^2\\
&=\|(r_{k+1}-u_{k+1})-(r_k-u_k)\|^2_{M+M^T-N^THN}\\
&=\|(r_{k+1}-u_{k+1})-(r_k-u_k)\|^2_{G}.
\end{aligned}
\end{equation*}
If $G$ is positive definite, we have
\begin{equation}
\label{eq5:19}
\|N(u_{k+1}-r_{k+1})\|_H^2<\|N(u_{k}-r_{k})\|_H^2. 
\end{equation}
Secondly, we turn to derive the convergence rate for Algorithm 1. Recalling Theorem \ref{th4} and according to the equivalence of different norms, there exists a constant $c>0$ such that
$$
\|u_{k+1}-u^*\|_H^2\leq 
\|u_{k}-u^*\|_H^2-c\|N(u_k-\tilde{u}_k)\|_H^2.
$$
Summing it over $k=0,1,\cdots,n$, we have
$$
\sum_{k=0}^n c\|N(u_k-r_k)\|_H^2\leq \|u_0-u^*\|_H^2.
$$
It follows from (\ref{eq5:19}) that $\{\|N(u_k-r_k)\|_H^2\}$ is monotonically nonincreasing. Therefore,
$$
(n+1)c\|N(u_n-r_n)\|_H^2\leq
\sum_{k=0}^n c\|N(u_k-r_k)\|_H^2\leq \|u_0-u^*\|_H^2,
$$
implying that
$$
\|N(u_n-r_n)\|_H^2\leq\frac{1}{(n+1)c} \|u_0-u^*\|_H^2,\forall u^*\in \operatorname{zer} P.
$$
This indicates the $O(\frac{1}{n})$ convergence rate for Algorithm \ref{alg:1}. The assertion of this theorem is proved.

\end{appendices}

\normalem    
\footnotesize 
\bibliographystyle{unsrt}
\bibliography{ref}
\end{document}